# Forward stagewise regression and the monotone lasso

**Trevor Hastie**[†]

*Departments of Statistics, and Health, Research & Policy
Sequoia Hall, Stanford University, CA 94305.
e-mail:* `hastie@stanford.edu`

**Jonathan Taylor**

*Department of Statistics
Sequoia Hall, Stanford University, CA 94305.
e-mail:* `jtaylor@stanford.edu`

**Robert Tibshirani**[‡]

*Departments of Health, Research & Policy, and Statistics
Sequoia Hall, Stanford University, CA 94305.
e-mail:* `tibs@stanford.edu`

**Guenther Walther**[§]

*Department of Statistics
Sequoia Hall, Stanford University, CA 94305.
e-mail:* `walther@stanford.edu`

**Abstract:** We consider the least angle regression and forward stagewise algorithms for solving penalized least squares regression problems. In Efron, Hastie, Johnstone & Tibshirani (2004) it is proved that the least angle regression algorithm, with a small modification, solves the lasso regression problem. Here we give an analogous result for incremental forward stagewise regression, showing that it solves a version of the lasso problem that enforces monotonicity. One consequence of this is as follows: while lasso makes optimal progress in terms of reducing the residual sum-of-squares per unit increase in $L_1$-norm of the coefficient $\beta$, forward stage-wise is optimal per unit $L_1$ arc-length traveled along the coefficient path. We also study a condition under which the coefficient paths of the lasso are monotone, and hence the different algorithms coincide. Finally, we compare the lasso and forward stagewise procedures in a simulation study involving a large number of correlated predictors.

[∗]The authors thank Steven Boyd, Jerome Friedman, Saharon Rosset, Ben van Roy and Ji Zhu for helpful discussions.

[†] Hastie was partially supported by grants DMS-0204612 and DMS-0505676 from the National Science Foundation, and grant 2R01 CA 72028-07 from the National Institutes of Health.

[‡]Tibshirani was partially supported by National Science Foundation Grant DMS-9971405 and National Institutes of Health Contract N01-HV-28183.

[§]Walther was partially supported by National Science Foundation Grant DMS-0505682 National Institutes of Health grant 5R33HL068522.

**Contents**

## 1. Introduction

The lasso (Tibshirani 1996) is a method for regularizing a least squares regression. Suppose we have predictor measurements $x_{ij}$, $j = 1, 2, \ldots, p$ and an outcome measurement $y_i$, observed for cases $i = 1, 2, \ldots N$. The lasso fits a linear model

$$f(x) = \beta_0 + \sum_{j=1}^{p} x_j \beta_j \qquad (1)$$

by solving the optimization problem

$$\min_\beta \sum_{i=1}^{N} \left( y_i - \beta_0 - \sum_{j=1}^{p} x_{ij} \beta_j \right)^2 \text{ subject to } \sum_{j=1}^{p} |\beta_j| \leq s \qquad (2)$$

If the tuning parameter $s \geq 0$ is large enough, this gives the ordinary least squares estimates. However, smaller values of $s$ produce shrunken estimates $\hat\beta$, often with many components equal to zero. Choosing $s$ can be thought of as choosing the number of predictors to include in a regression model. Thus the lasso can select predictors like subset selection methods. However, since it is a smooth optimization problem, it is less variable than subset selection, and can be applied to much larger problems (large in $p$). Chen, Donoho & Saunders (1998) developed related technology in the context of signal processing.

The criterion (2) leads to a quadratic programming problem for each $s$, and thus standard numerical-analysis methods can be used to solve it. Figure 1 shows an example, based on simulated data with 10 predictors (details of the data generation and model are given later). The top panel shows the coefficient profiles of the lasso solutions, as the bound $s = \sum |\beta_j|$ is increased from 0 up to the point where the full least squares solutions are obtained (right end of figure).

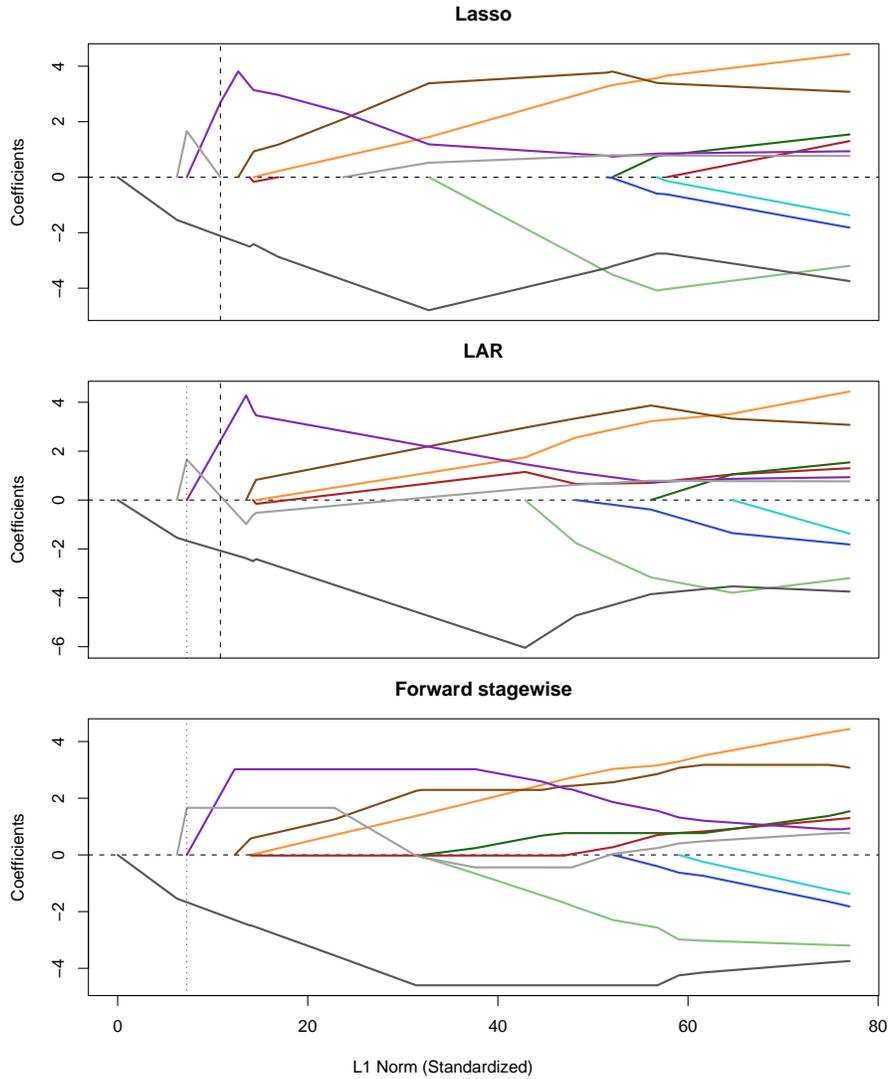

FIG 1. *Coefficient profiles, simulated example. The $L_1$ norm is computed on the standardized variables. The coefficients are given on their original scale, on which the details are more visible. The lasso starts to differ from LAR at the broken vertical line, when the gray coefficient passes through zero. Forward stagewise starts to differ from LAR and lasso at the dotted vertical line, where the gray coefficient goes flat instead of turning back towards zero.*

Notice the piecewise linear nature of the lasso profiles. Efron et al. (2004) exploited this fact to derive a simple algorithm — *least angle regression* — for simultaneously solving the entire set of lasso problems (all values of $s$). This work was motivated by an observation in Hastie, Tibshirani & Friedman (2001, Section 10.12.2) that the lasso profile bore a striking similarity to the coefficient profile produced by a version of boosting for linear models, which they named the *incremental forward stagewise* algorithm (hereafter $FS_\epsilon$).

This $FS_\epsilon$ algorithm (see Algorithm 1) creates a coefficient profile as follows: at each step it increments the coefficient of that variable most correlated with the current residuals by an amount $\pm\epsilon$, with the sign determined by the sign of the correlation. Efron et al. (2004) in fact considered the limiting version of this algorithm, with $\epsilon \downarrow 0$, which also has piecewise linear coefficient paths. We refer to this as the *infinitesimal forward stagewise* algorithm, hereafter $FS_0$ or simply *forward stagewise*. Efron et al. (2004) showed that under certain conditions, this $FS_0$ path is identical to the lasso path. However, for most problems they are different (e.g. Figure 1), and sometimes strikingly so (Figure 7). The $FS_0$ paths are much smoother than the lasso paths.

The primary result in this paper is the characterization of $FS_0$ as a monotone version of the lasso, in a sense to be described in Section 3. As such, it is a more restricted version of the lasso, and hence the additional smoothness. Because of the monotonicity, the criterion cannot be defined pointwise (as lasso can), but instead defines the entire path via a differential equation.

In Section 4 we generalize this characterization to loss functions other than squared error. Section 5 considers other candidate criteria, and Section 6 examines conditions under which the lasso and $FS_0$ are the same. Section 7 compares the two procedures in a simulation study. Some proofs are given in the Appendix.

## 2. Background: The LARS Algorithm

Hastie et al. (2001) showed that the solution path for the lasso is strikingly similar to that of a simplified version of "boosting". Boosting is an adaptive, non-linear, function-fitting method that has received much attention in the past ten years (Schapire & Freund 1997, Schapire, Freund, Bartlett & Lee 1998, Friedman, Hastie & Tibshirani 2000). In modern versions of boosting, the set of "variables" is a large space of binary trees, which are selected, shrunk, and added to the current model. In their simplification, Hastie et al. (2001) replaced boosted trees by an incremental forward stagewise algorithm for linear regression, reproduced here as Algorithm 1.

The Least Angle Regression algorithm (LAR, see Algorithm 2) of Efron et al. (2004) was initially intended as the limiting version $FS_0$ of $FS_\epsilon$. As we explain below, LAR is different from both $FS_0$ and lasso, but both can be obtained through simple modifications. The bottom panel of Figure 1 shows the coefficient profiles of $FS_0$. Notice that they are similar to lasso and LAR, but tend to be smoother.[1]

**Algorithm 1** *Incremental Forward Stagewise Regression: $FS_\epsilon$*
1. Start with $\mathbf{r} = \mathbf{y} - \bar{\mathbf{y}}$, $\beta_1, \beta_2, \ldots \beta_p = 0$.
2. Find the predictor $\mathbf{x}_j$ most correlated with $\mathbf{r}$.
3. Update $\beta_j \leftarrow \beta_j + \delta_j$, where $\delta_j = \epsilon \cdot \text{sign}[\text{corr}(\mathbf{r}, \mathbf{x}_j)]$;
4. Update $\mathbf{r} \leftarrow \mathbf{r} - \delta_j \mathbf{x}_j$, and repeat steps 2 and 3 until no predictor has any correlation with $\mathbf{r}$.

LAR is a kind of "democratic" alternative to a version of the commonly used forward-stepwise regression algorithm. Forward-stepwise regression starts with all coefficients equal to zero, and then builds a sequence of models by successively including one variable at a time, and updating the least-squares fit. The version we consider here enters at each stage the variable most correlated with the residuals.[2] This process is repeated until all $p$ predictors have been entered, or the residuals are zero.

LAR uses a similar strategy, but only enters "as much" of a predictor as it "deserves": the coefficient of the predictor is increased only up to the point where some other predictor has as much correlation with the current residual. This new predictor is entered, and the process is continued. Algorithm 2 gives more details.

**Algorithm 2** *LARS: Least Angle Regression*
1. Standardize the predictors to have mean zero and variance 1. Start with the residual $\mathbf{r} = \mathbf{y} - \bar{\mathbf{y}}$, $\beta_1, \beta_2, \ldots \beta_p = 0$.
2. Find the predictor $\mathbf{x}_j$ most correlated with $\mathbf{r}$.
3. Move $\beta_j$ from 0 towards its least-squares coefficient $\langle \mathbf{x}_j, \mathbf{r} \rangle$, until some other competitor $\mathbf{x}_k$ has as much correlation with the current residual as does $\mathbf{x}_j$.
4. Move $(\beta_j, \beta_k)$ in the direction defined by their joint least squares coefficient of the current residual on $(\mathbf{x}_j, \mathbf{x}_k)$, until some other competitor $\mathbf{x}_l$ has as much correlation with the current residual.
5. Continue in this way until all $p$ predictors have been entered. After $p$ steps, we arrive at the full least-squares solution.

The profiles for LAR are shown in the middle panel of Figure 1. They look similar to the lasso solutions, especially in the beginning. The first discrepancy is at the place marked by a vertical broken line in the lasso profiles. The LAR

---
[1] The name "LARS", derived from least angle regression and lasso, is the name we have given to our algorithm that implements LAR, lasso and $FS_0$.

[2] This can differ from the more traditional version, which includes the variable that leads to the largest drop in residual sum-of-squares.

profile passes through zero at this point, while the lasso profile hits zero, and stays there. This similarity is no coincidence. It turns out that with one modification, the LAR procedure exactly produces the set of lasso solutions for all $s$. The modification needed is as follows:

---
**Algorithm 2a** *LARS: lasso Modification.*

5a. If a non-zero coefficient hits zero, drop it from the active set and recompute the current joint least squares direction.

---

The LARS(lasso) algorithm is extremely efficient, requiring the same order of computation as that of a single least squares fit using the $p$ predictors. Least angle regression always takes $p$ steps to get to the full least squares estimates. The lasso path can have more than $p$ steps, although the two are often quite similar. Algorithm 2a is an efficient way of computing the solution to any lasso problem, especially when $N \ll p$ (Donoho & Tsaig 2006). Osborne, Presnell & Turlach (2000) also discovered a piecewise-linear path for computing the lasso, which they called a *homotopy* algorithm.

Efron et al. (2004) showed that another variant of the LAR algorithm gives $FS_0$; see also Theorem 1 on page 10. Suppose we have reached a point when a variable enters the active set $\mathcal{A}$; all variables $\mathbf{x}_j$, $j \in \mathcal{A}$ have correlation equal in magnitude with the current residual $\mathbf{r}$. Then the new LAR direction is defined by the least-squares fit of $\mathbf{r}$ on $\mathbf{X}_\mathcal{A}$. The modification needed to achieve $FS_0$ replaces this by a type of *non-negative* least squares direction.

---
**Algorithm 2b** *LARS: $FS_0$ Modification.*

4. Find the new direction by solving the constrained least squares problem

$$\min_b ||\mathbf{r} - \mathbf{X}_\mathcal{A} b||_2^2 \text{ subject to } b_j s_j \geq 0, \ j \in \mathcal{A},$$

where $s_j$ is the sign of $\langle \mathbf{x}_j, \mathbf{r} \rangle$.

---

The constraints arise from the fact that in step 3 of the incremental forward stagewise procedure, the coefficient of each predictor is increased in the direction of its correlation with the current residual. In Figure 1 LAR and $FS_0$ start to differ (vertical dotted line) when the third variable enters $\mathcal{A}$.

The top left panel of Figure 3 on page 12 shows the residual sum of squares (RSS) for each of the three procedures, as a function of the $L_1$ norm of the coefficient vector. As expected, the lasso curve lies below the other two, because it decreases the RSS the fastest per unit increase in $L_1$ norm. The right panel plots RSS against the $L_1$ arc-length of the coefficient profile; here $FS_0$ wins — a point we enlarge on in the next section.

In summary, we see that the forward stagewise and LAR algorithms "nearly" solve the $L_1$-penalized regression problem. It is natural to ask: what problems are the forward stagewise and LAR algorithms solving? Keith Knight asked that question in the discussion of Efron et al. (2004).

We provide some answers to the first question in this paper. We make three main contributions:

- we characterize forward stagewise as a monotone version of the lasso, in an extended space of variables consisting of each variable and its negative.
- we study a condition under which the profiles of all three methods are monotone, and hence the three methods coincide.
- we compare the lasso and forward stagewise procedures in a simulation study involving a large number of correlated predictors.

## 3. Forward Stagewise and the Monotone Lasso

In this section we consider an expanded representation of the lasso problem which facilitates a clearer understanding of the forward stagewise procedure. For each predictor $x_j$, we include its *negative* version $-x_j$, resulting in an expanded data set with $2p$ predictors. In matrix notation we create an expanded data matrix $\tilde{\mathbf{X}} = [\mathbf{X} : -\mathbf{X}]$. In this framework the lasso problem becomes

$$= \min_{\beta_0, \beta_j^+, \beta_j^-} \sum_{i=1}^{n} \left( y_i - \beta_0 - \left[ \sum_{j=1}^{p} x_{ij} \beta_j^+ - \sum_{j=1}^{p} x_{ij} \beta_j^- \right] \right)^2 \quad (3)$$

$$\text{subject to } \beta_j^+, \beta_j^- \geq 0 \; \forall j \text{ and } \sum_{j=1}^{p} (\beta_j^+ + \beta_j^-) \leq s$$

In what follows we will sometimes suppress the constant term $\beta_0$, which can always be removed once and for all by centering all the variables. Since each value of the bound $s$ characterizes a solution, we can use $s$ to index the solution $\beta(s)$; whenever the constraint is active, the solution satisfies $||\beta(s)||_1 = s$, and we say the solution profile is parametrized by $L_1$-norm. Problem (3) is equivalent to a standard representation for solving the lasso problem by quadratic programming, and the KKT conditions ensure that at most one of $\hat{\beta}_j^+$ and $\hat{\beta}_j^-$ are greater than zero at the same time. Hence by augmenting the data with the negative of the variables, the *positive* lasso in the enlarged space is equivalent to the original lasso problem. Figure 2 (top pair of panels) shows the coefficient paths of the positive and negative variables for the lasso solution in Figure 1.

In the lower pair of plots, an additional constraint is imposed on this sequence of lasso problems: the coefficient paths are constrained to be *monotone nondecreasing*. These monotone paths are exactly equivalent to the paths of the forward-stagewise algorithm. By this we mean that the *collapsed* versions of the paths (subtracting the coefficients for the negative versions of the variables, from the corresponding coefficients for the positive versions) are exactly the forward-stagewise paths in the lower panel in Figure 1.

This leads us to characterize the forward-stagewise algorithm as a monotone version of the lasso. These extra restrictions are an additional form of regularization, leading to smoother coefficient profiles.

This expanded space of variables creates a more natural analog of boosting, which operates in a large dictionary of binary trees. For every tree, its negative is also available.[3] In the expanded space the equivalent of Algorithm 1 is given in Algorithm 3. It is obvious that Algorithm 3 generates monotone coefficient

---

**Algorithm 3** Monotone Incremental Forward Stagewise Regression

1. Start with $\mathbf{r} = \mathbf{y} - \bar{\mathbf{y}}$, $\beta_1, \beta_2, \ldots \beta_{2p} = 0$.
2. Find the predictor $\tilde{\mathbf{x}}_j$ most positively correlated with $\mathbf{r}$.
3. Update $\beta_j \leftarrow \beta_j + \epsilon$.
4. Update $\mathbf{r} \leftarrow \mathbf{r} - \epsilon \tilde{\mathbf{x}}_j$, and repeat steps 2 and 3 until no predictor has any correlation with $\mathbf{r}$.

---

paths, indexed by the number of steps $m$, or the total distance stepped $t = m \cdot \epsilon$. Drawing on the results of Efron et al. (2004), we show in Theorem 1 that the limit as $\epsilon \downarrow 0$ leads exactly to the monotone representation as in Figure 2. First we define the notion of $L_1$ arc-length.

**Definition 1.** *Suppose $\beta(t)$ is a one-dimensional differentiable curve in $t \geq 0$, with $\beta(0) = 0$. The $L_1$ arc-length of $\beta(t)$ in $[0, t]$ is given by*

$$\text{TV}(\beta, t) = \int_0^t ||\dot{\beta}(s)||_1 ds, \qquad (4)$$

*where $\dot{\beta}(s) = \partial \beta(s)/\partial s$.*

We have named the arc-length "TV" for *total-variation*; the $L_1$ arc-length of $\beta(t)$ up to time $t$ is the sum of total variation measures for each of its coordinate functions, and is a measure of roughness of the curve.

For a piecewise-differentiable continuous curve, the arc-length is the sum of the arc-lengths of the differentiable pieces. The following lemma is easily proved:

**Lemma 1.** *If the coordinates of $\beta(t)$ are monotone and piecewise differentiable in $t$, then $\text{TV}(\beta, t) = ||\beta(t)||_1$.*

Hence the arc-length and $L_1$ norm for a monotone coefficient profile are the same.

Although it is convenient to use this expanded representation $\tilde{\mathbf{X}}$, we can always collapse to the original representation $\mathbf{X}$. The coefficients in the original representation are simply the paired differences $\beta_j(t) = \beta_j^+(t) - \beta_j^-(t)$. Note that

- The $L_1$ norm for the lasso coefficients is the same in either representation, since only one coefficient in each pair is non-zero at a time (see the proof for part 1 of Theorem 1).
- The $L_1$ arc-length for $\text{FS}_0$ is the $L_1$ norm in the expanded representation, and is equal to the $L_1$ arc-length in the original representation. This is NOT the same as the $L_1$ norm in the original representation.

---

[3]We can think of a tree as a variable; the $N$ values are obtained by passing the training data down to the terminal nodes.

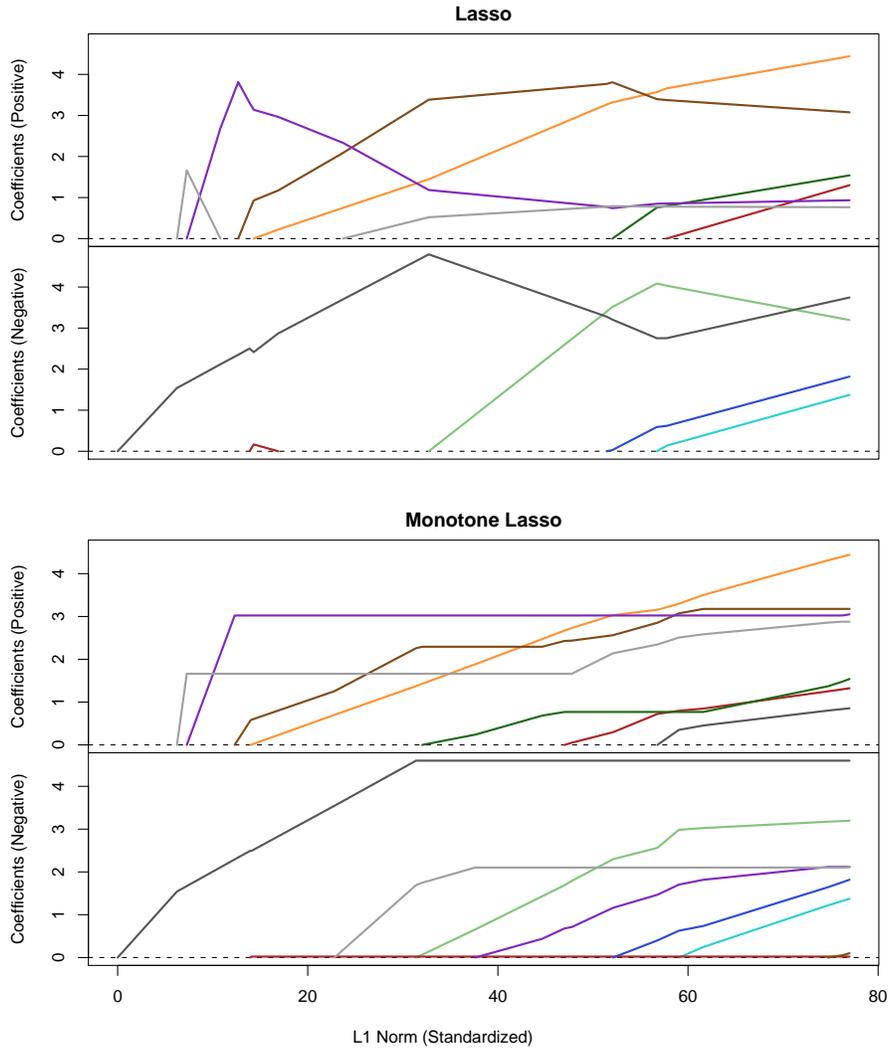

FIG 2. *Expanded coefficient profiles, simulated example. The $L_1$ norm is computed on the standardized variables. The coefficients are given on their original scale, on which the details are more visible.*

Every point along the lasso path is the solution to a convex optimization problem. Unfortunately, the monotonicity restriction of the forward stagewise path appears to preclude such a succinct characterization. Alternatively, we can show that the lasso path is the solution to a differential equation, which characterizes the path in terms of a series of optimal moves. We then show that the forward stagewise path is the solution to a closely related differential equation, which restricts these optimal moves to be monotone. In the remainder of this section:

- we characterize the forward stagewise path in terms of a sequence of monotone *moves*, and compare these moves to the less restrictive moves of the lasso (Theorem 1);
- this leads us to define the *monotone lasso* — a path defined by a differential equation — with the derivatives giving the move directions from the current position (Definitions 2–3). The lasso can also be characterized as a solution to a related differential equation;
- we show that the monotone lasso is locally optimal in terms of *arc-length*— it makes the optimal move per unit increase in arc-length of the coefficient profile. The lasso makes the optimal move per unit increase in the $L_1$ norm of the coefficients (Theorem 2);
- we show that the forward stagewise algorithm computes the solution to the monotone lasso criterion. (Proposition 1).

We then generalize these results for other loss functions in Section 4.

**Theorem 1.** *Let $\beta^0 \in \mathbb{R}^{2p}$ be a point either on the lasso or forward stagewise path in the expanded-variable space, and let $\mathcal{A}$ be the* active *set of variables achieving the maximal correlation with the current residual $\mathbf{r} = \mathbf{y} - \tilde{\mathbf{X}}\beta^0$.*

1. *The lasso coefficients move in a direction given by the coefficients of the least squares fit of $\tilde{\mathbf{X}}_\mathcal{A}$ on $\mathbf{r}$.*
2. *The forward stagewise coefficients move in a direction given by the coefficients of the* non-negative *least squares fit of $\tilde{\mathbf{X}}_\mathcal{A}$ on $\mathbf{r}$.*

*In either case only the coefficients in $\mathcal{A}$ change, and this fixed direction is pursued until the first of the following events occurs:*

(a) *a variable not in $\mathcal{A}$ attains the maximal correlation and joins $\mathcal{A}$;*
(b) *The coefficient of a variable in the active set reaches 0, at which point it leaves $\mathcal{A}$ (lasso only);*
(c) *the residuals match those of the unrestricted least squares fit.*

*When (a) or (b) occur, the direction is recomputed.*

The proof of this theorem can be assembled from the results proved in Efron et al. (2004). For convenience we give a simple proof in the appendix, using convex optimality conditions (see also Rosset & Zhu (2004)).

Theorem 1 leads us to define the *monotone lasso* as the solution to a differential equation, which is characterized in terms of its positive path derivatives. The FS$_0$ algorithm computes this solution.

Theorem 1 is stated in terms of a point $\beta^0$ on the lasso/FS$_0$ paths. In fact these moves can be defined starting from any value $\beta^0$.

**Definition 2.** *Let $\beta \in \mathbb{R}^{2p}$ be any coefficient for a linear model in the expanded variable set, and let $\mathbf{r} = \mathbf{y} - \tilde{\mathbf{X}}\beta$. Let $\mathcal{A}$ be the active set of variables achieving maximal correlation with $\mathbf{r}$.*

1. *The* **lasso move direction** *$\rho_l(\beta) : \mathbb{R}^{2p} \mapsto \mathbb{R}^{2p}$ is defined*

$$\rho_l(\beta) = \begin{cases} 0 & \text{if } \tilde{\mathbf{X}}^T\mathbf{r} = 0 \\ \theta/\sum_j \theta_j & \text{otherwise,} \end{cases} \quad (5)$$

   *with $\theta_j = 0$ except for $j \in \mathcal{A}$, where $\theta_\mathcal{A}$ is the least squares coefficient of $\mathbf{r}$ on $\tilde{\mathbf{X}}_\mathcal{A}$.*

2. *The* **monotone lasso move direction** *$\rho_{ml}(\beta) : \mathbb{R}^{2p} \mapsto \mathbb{R}^{2p}$ is defined*

$$\rho_{ml}(\beta) = \begin{cases} 0 & \text{if } \tilde{\mathbf{X}}^T\mathbf{r} = 0 \\ \theta/\sum_j \theta_j & \text{otherwise,} \end{cases} \quad (6)$$

   *with $\theta_j = 0$ except for $j \in \mathcal{A}$, where $\theta_\mathcal{A}$ is the non-negative least squares coefficient of $\mathbf{r}$ on $\tilde{\mathbf{X}}_\mathcal{A}$.*

The normalizations in (5) and (6) are not essential, but turn out to be convenient when we parametrize the coefficient paths later in this section.

Figure 3 shows the residual-sum-of-squares (RSS) curves for the lasso and forward stagewise algorithms, applied to our simulation example. It appears in this example that lasso decreases RSS most rapidly as a function of the $L_1$ norm of the coefficients $||\beta(t)||_1$, while forward stagewise wins in terms of $L_1$ arc-length. It turns out that this is always the case, and is a characterization of the local optimality for each of the procedures.

**Theorem 2.** *Let $\beta^0 \in \mathbb{R}^{2p}$ be a coefficient vector in the expanded-variable space. Then the lasso/monotone lasso move directions defined in Definition 2 are optimal in the sense that*

1. *A lasso move decreases the residual sum of squares at the optimal quadratic rate with respect to the $L_1$ coefficient norm;*
2. *A monotone-lasso move decreases the residual sum of squares at the optimal quadratic rate with respect to the coefficient $L_1$ arc-length.*

There is some intuition in this distinction when we think of forward stagewise as a form of boosting. There we pay a cost in terms of effort for any move we make (number of trees), which is captured by arc-length. With the lasso we get rewarded for decreasing a coefficient towards zero. The monotonicity constraint also results in much smoother coefficient profiles, and hence shorter arc-lengths. Zhao & Yu (2004) propose a modification to boosting to allow the backtracking needed to make FS$_0$ coincide with lasso.

Our proof follows closely the material in Section 6 of Efron et al. (2004). Since the directions are fixed while $\mathcal{A}$ is fixed, the paths are piecewise linear, and hence the residual-sum-of-squares curves are piecewise quadratic.

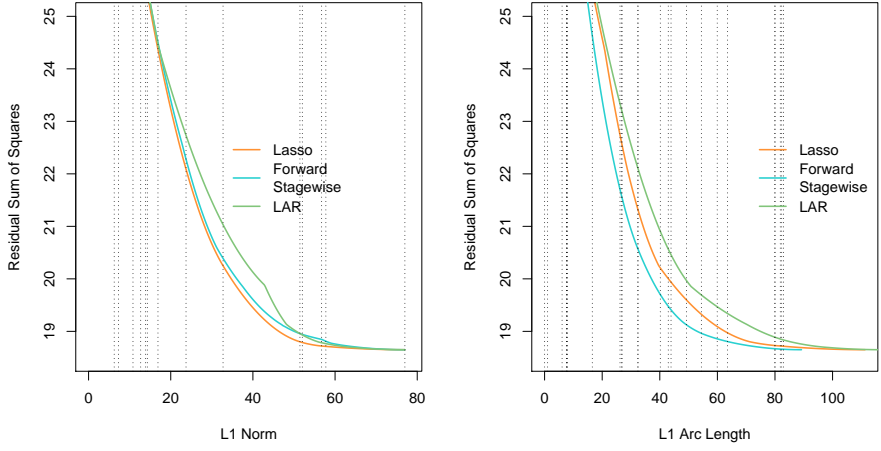

FIG 3. *The RSS for our simulation example, as a function of the $L_1$ norm (left panel) and arc length (right panel) of the coefficient paths for lasso, forward stagewise, and Least Angle Regression.*

*Proof of Theorem 2: lasso.* Consider a move in direction $d$ from $\beta^0 : \beta^0 + \gamma \cdot d$. Define

$$R(\gamma) = ||\mathbf{y} - \tilde{\mathbf{X}}(\beta^0 + \gamma \cdot d)||_2^2; \qquad (7)$$
$$T(\gamma) = (\beta^0 + \gamma \cdot d)^T \mathbf{1}, \qquad (8)$$

where $\mathbf{1} = (1, \ldots, 1)^T$. Assuming $d_j > 0$ when $\beta_j^0 = 0$, $T(\gamma)$ is the $L_1$ norm of the changed coefficient. We then compute the path derivative

$$U(\gamma) = \frac{\partial R}{\partial T} = \frac{\partial R(\gamma)}{\partial \gamma} \Big/ \frac{\partial T(\gamma)}{\partial \gamma} \qquad (9)$$

$$= -2\frac{d^T \tilde{\mathbf{X}}^T \left(\mathbf{y} - \tilde{\mathbf{X}}(\beta^0 + \gamma \cdot d)\right)}{d^T \mathbf{1}}; \qquad (10)$$

$$U(\gamma)|_{\gamma=0} = -2\frac{d^T \tilde{\mathbf{X}}^T \mathbf{r}}{d^T \mathbf{1}}, \qquad (11)$$

where $\mathbf{r} = \mathbf{y} - \tilde{\mathbf{X}}\beta^0$. Since $\tilde{\mathbf{x}}_j^T \mathbf{r} = C$ for all $j \in \mathcal{A}$, the maximal-correlation active set, this derivative is minimized by allowing only those elements $d_j$ with $j \in \mathcal{A}$ to be nonzero. For any such $d = d_\mathcal{A}$, the derivative is $U(0) = -2C$. Among these, we seek the $d_\mathcal{A}$ with smallest Hessian.

$$\frac{\partial^2 R}{\partial T^2} = \frac{\partial U}{\partial T} = \frac{\partial U}{\partial \gamma} \Big/ \frac{\partial T}{\partial \gamma}$$
$$= 2\frac{d_\mathcal{A}^T \tilde{\mathbf{X}}_\mathcal{A}^T \tilde{\mathbf{X}}_\mathcal{A} d_\mathcal{A}}{(d_\mathcal{A}^T \mathbf{1})^2} \qquad (12)$$

Minimizing this Rayleigh-quotient is equivalent to minimizing

$$d_{\mathcal{A}}^T \tilde{\mathbf{X}}_{\mathcal{A}}^T \tilde{\mathbf{X}}_{\mathcal{A}} d_{\mathcal{A}} \quad \text{subject to } d_{\mathcal{A}}^T \mathbf{1} = 1.$$

It is straightforward to show that the solution is given by $d_{\mathcal{A}} \propto (\tilde{\mathbf{X}}_{\mathcal{A}}^T \tilde{\mathbf{X}}_{\mathcal{A}})^{-1} \mathbf{1}$. But since $\tilde{\mathbf{X}}_{\mathcal{A}}^t \mathbf{r} = C \cdot \mathbf{1}$, this is equivalent to the lasso move.

Hence the sequence of lasso moves result in an optimal piecewise-quadratic RSS drop-off curve as a function of $L_1$ norm.

*Proof of Theorem 2: monotone lasso.*

The increment in the $L_1$-arc-length of the path $\beta^0 + \gamma \cdot d$ (starting from $\beta^0$) is easily seen to be

$$L(\gamma) = \gamma ||d||_1. \tag{13}$$

Similar to (9), we get

$$V(\gamma) = \frac{\partial R}{\partial L} = \frac{\partial R(\gamma)}{\partial \gamma} \bigg/ \frac{\partial L(\gamma)}{\partial \gamma} \tag{14}$$

$$= -2 \frac{d^T \tilde{\mathbf{X}}^T (\mathbf{y} - \tilde{\mathbf{X}}(\beta + \gamma \cdot d))}{||d||_1}. \tag{15}$$

$$V(\gamma)|_{\gamma=0} = -2 \frac{d^T \tilde{\mathbf{X}}^T \mathbf{y}}{||d||_1}. \tag{16}$$

This is minimized by selecting $d_j \geq 0$ for $j \in \mathcal{A}$, again with minimizing value $-2C$. The Hessian is

$$\frac{\partial^2 R}{\partial L^2} = \frac{\partial V}{\partial \gamma} \bigg/ \frac{\partial L}{\partial \gamma}$$

$$= 2 \frac{d_{\mathcal{A}}^T \tilde{\mathbf{X}}_{\mathcal{A}}^T \tilde{\mathbf{X}}_{\mathcal{A}} d_{\mathcal{A}}}{||d||_1^2}, \tag{17}$$

which we would like to minimize subject to $d_j \geq 0$. This is equivalent to the optimization problem

$$\min_d d^T \tilde{\mathbf{X}}_{\mathcal{A}}^T \tilde{\mathbf{X}}_{\mathcal{A}} d \quad \text{subject to } d_j \geq 0, \sum_{j \in \mathcal{A}} d_j = 1. \tag{18}$$

It is straightforward to show via the KKT conditions for this quadratic programming problem that the solution is identical to the solution for $\rho$ in (53) in the appendix, which is the direction given by a non-negative least-squares fit of $\mathbf{r}$ to $\tilde{\mathbf{X}}_{\mathcal{A}}$ (the forward-stagewise move). □

The graphs in Figure 3 suggest that the gap is bigger as a function of arc-length than norm. This is in fact the case, as can be seen in the proof of Theorem 2. As a function of norm, starting from the same point, the downward gradient is the same for both lasso and $FS_0$, but the Hessian is smaller for lasso. As a function of arc-length, the gradient for lasso can be larger than for $FS_0$, if some of the $d_j$ are negative.

Armed with the lasso and monotone lasso move directions from definition 2, we can now characterize paths as solutions to differential equations.

**Definition 3.** *The monotone lasso coefficient path $\beta(\ell)$ for a dataset $\tilde{\mathbf{X}} = \{\mathbf{X}, -\mathbf{X}\}$ is the solution to the differential equation*

$$\frac{\partial \beta}{\partial \ell} = \rho_{ml}(\beta(\ell)), \tag{19}$$

*with initial condition $\beta(0) = 0$. Since (19) is piecewise continuous, this path is continuous and piecewise differentiable.*

Because the directions $\rho_{ml}(\beta)$ defined in (6) are standardized to have unit $L_1$ norm, the solution curve is *unit $L_1$-speed*, and hence is parametrized by $L_1$ arc length.

In order to solve (19), we need to track the entire path; this solution is provided by the forward-stagewise algorithm.

**Proposition 1.** *The forward-stagewise algorithm for a dataset $\tilde{\mathbf{X}} = \{\mathbf{X}, -\mathbf{X}\}$ and square-error loss computes the monotone lasso path $\beta(\ell)$; it starts at 0, and then increments the coefficients continuously according to the monotone lasso moves (6). Specifically*

**Initialize** Set $\beta(0) = 0$, $\ell_0 = 0$, and $\rho_0 = \rho_{ml}(0)$, with corresponding active set $\mathcal{A}_0$.

**For** $j = 0, 1, 2, \ldots$

1. Let $\beta(\ell) = \beta(\ell_j) + (\ell - \ell_j) \cdot \rho_j$, $\ell \in [\ell_j, \ell_{j+1}]$, where $\ell_{j+1}$ is the value of $\ell > \ell_j$ at which $\mathcal{A}_j$ changes to $\mathcal{A}_{j+1}$.
2. Compute $\rho_{j+1} = \rho_{ml}(\beta(\ell_{j+1}))$.
3. If $\rho_{j+1} = 0$ exit, and $\beta(\ell)$ is defined on $[0, L]$, with $L = \ell_{j+1}$.

Proposition 1 follows from Theorem 1. We can characterize the lasso path in a similar fashion.[4]

**Proposition 2.** *The lasso coefficient path $\beta(\ell)$ for a dataset $\tilde{\mathbf{X}} = \{\mathbf{X}, -\mathbf{X}\}$ is the solution to the differential equation*

$$\frac{\partial \beta}{\partial \ell} = \rho_l(\beta(\ell)), \tag{20}$$

*with initial condition $\beta(0) = 0$. Since (20) is piecewise continuous, this path is continuous and piecewise differentiable.*

The normalization of $\rho_l$ defined in (5) guarantees that the solution path is parametrized by $L_1$ norm (since the coefficients are non-negative).

The characterizations above draw on the similarities between the lasso and monotone lasso. The characterization of the monotone lasso falls slightly short of that of the lasso for the following reasons.

---

[4]Since the lasso path is defined in Tibshirani (1996) as the solution to a convex optimization problem, this alternative characterization is a proposition (unlike Definition 3), and follows from Theorem 1.

- We can define a lasso solution explicitly at any given point on the path, as the solution to an optimization problem (2); we are unable to do this for the monotone lasso.
- When $p < n$, both the lasso and monotone-lasso paths end in the unrestricted least-squares solution. When $p > n$, any least squares solution has zero residuals, with infinitely many solution coefficients $\beta$. The lasso path leads to the unique zero-residual solution having minimum $L_1$ norm. By construction the monotone lasso path also produces a unique zero-residual solution in these circumstances, but we are unable to characterize it further.

## 4. Forward Stagewise for General Convex Loss Functions

Gradient boosting (Friedman 2001, Hastie et al. 2001) is often used with loss functions other than squared error; typical candidates are the binomial log-likelihood or the "Adaboost" loss for binary classification problems. Our linear-model simplification is also applicable there. As a concrete example, consider the linear logistic regression model in the expanded space:

$$\log \frac{\Pr(y=1|\tilde{x})}{\Pr(y=0|\tilde{x})} = \tilde{x}^T \beta \qquad (21)$$

$$= \eta(\tilde{x}). \qquad (22)$$

The negative of the binomial log-likelihood is given by

$$L(\beta) = -\sum_{i=1}^{n} \left[ y_i \log p_i + (1 - y_i) \log(1 - p_i) \right], \qquad (23)$$

where

$$p_i = \frac{e^{\tilde{x}_i^T \beta}}{1 + e^{\tilde{x}_i^T \beta}}. \qquad (24)$$

More generally, consider the case where we have a linear model $\eta(\tilde{x}) = \tilde{x}^T \beta$, and a loss function of the form

$$L(\beta) = \sum_{i=1}^{n} l(y_i, \eta(\tilde{x}^T \beta)). \qquad (25)$$

The analog of Algorithm 3 for this general case is given in Algorithm 4.

For the binomial case, the negative gradient in step 2. is $-\partial L/\partial \beta_j = \tilde{\mathbf{x}}_j^T(\mathbf{y} - \mathbf{p})$, where $\mathbf{y}$ is the 0/1 response vector, and $\mathbf{p}$ the vector of fitted probabilities.

We can apply the same logic used in forward stagewise with squared error loss in this situation, by using a quadratic approximation to the loss at the current $\beta^0$:

$$L(\beta) \approx L(\beta^0) + \frac{\partial L}{\partial \beta}\Big|_0 (\beta - \beta^0) + \frac{1}{2}(\beta - \beta^0)^T \frac{\partial^2 L}{\partial \beta \partial \beta^T}\Big|_0 (\beta - \beta^0). \qquad (26)$$

**Algorithm 4** Generalized Monotone Incremental Forward Stagewise Regression
1. Start with $\beta_1, \beta_2, \ldots \beta_{2p} = 0$.
2. Find the predictor $x_j$ with largest negative gradient element $-\partial L/\partial \beta_j$, evaluated at the current predictor $\eta$.
3. Update $\beta_j \leftarrow \beta_j + \epsilon$.
4. Update the predictors $\eta(x_i) = \tilde{x}_i^T \beta$, and repeat steps 2 and 3 many times

The two derivatives in this case are

$$\frac{\partial L}{\partial \beta}\Big|_0 = \tilde{\mathbf{X}}^T \mathbf{u}^0 \tag{27}$$

$$\frac{\partial^2 L}{\partial \beta \partial \beta^T}\Big|_0 = \tilde{\mathbf{X}}^T \mathbf{W}^0 \tilde{\mathbf{X}}, \tag{28}$$

where

$$u_i^0 = \frac{\partial l(y_i, \eta)}{\partial \eta}\Big|_{\eta = \tilde{x}_i^T \beta^0}, \tag{29}$$

and the diagonal matrix $\mathbf{W}^0$ has entries

$$W_{ii}^0 = \frac{\partial^2 l(y_i, \eta)}{\partial \eta^2}\Big|_{\eta = \tilde{x}_i^T \beta^0}. \tag{30}$$

In the case of logistic regression $W_{ii}^0 = p_i^0(1 - p_i^0)$, where $p_i^0$ are the current probabilities, and $u_i^0 = -(y_i - p_i^0)$. Minimizing (26) gives the Newton update

$$\begin{aligned} \delta &= \beta - \beta^0 \\ &= -(\tilde{\mathbf{X}}^T \mathbf{W}^0 \tilde{\mathbf{X}})^{-1} \tilde{\mathbf{X}}^T \mathbf{u}^0, \end{aligned} \tag{31}$$

which can be expressed as the coefficients from a weighted least squares fit of $\tilde{\mathbf{X}}$ on $-\mathbf{W}^{0^{-1}} \mathbf{u}^0$, with weights $\mathbf{W}^0$.

**Definition 4.** *The monotone lasso move direction $\rho_{ml}(\beta, L)$ at a point $\beta^0$, with expanded data $\tilde{\mathbf{X}}$, and with loss function $L$ is:*

1. *Compute $\tilde{\mathbf{X}}^T \mathbf{u}^0$; if all elements are zero, return $\rho = 0$.*
2. *Establish the active set $\mathcal{A}$ of indices for which $-\tilde{\mathbf{x}}_j^T \mathbf{u}^0 = \max_{k=1}^{2p} -\tilde{\mathbf{x}}_k^T \mathbf{u}^0$.*
3. *Let $\hat{\delta}$ be the coefficients from a weighted, positive, least squares fit of $\tilde{\mathbf{X}}_{\mathcal{A}}$ on $-\mathbf{W}^{0^{-1}} \mathbf{u}^0$, with weights $\mathbf{W}^0$.*
4. *Define*

$$\rho_j = \begin{cases} \hat{\delta}_j / \sum_j \hat{\delta}_j & \text{if } j \in \mathcal{A} \\ 0 & \text{otherwise.} \end{cases} \tag{32}$$

It is easy to check that this definition coincides with the definition for squared-error loss. Unlike there, it will in general not be piecewise constant. We can expect it to be piecewise smooth, with breaks when the active sets change.

**Definition 5.** *The monotone lasso coefficient path $\beta(\ell)$ for a dataset $\tilde{\mathbf{X}} = \{\mathbf{X}, -\mathbf{X}\}$ and loss $L$ is the solution to the differential equation*

$$\frac{\partial \beta}{\partial \ell} = \rho_{ml}(\beta(\ell), L), \tag{33}$$

*with initial condition $\beta(0) = 0$.*

The definitions are exactly analogous for generalizations of the lasso.

Unlike for squared error loss, the solution paths are in general piecewise smooth (but nonlinear), and so efficient exact path algorithms are not available. Rosset (2005) show that as long as the loss function is quadratic, piecewise linear or a mixture of both, then the paths will be piecewise linear, and can be tracked.

For these general convex-loss lasso problems, lasso solutions are always available at any point along the path. Park & Hastie (2006) develop efficient algorithms for obtaining the lasso path for the *generalized linear model* family of loss functions (including logistic regression).

For the monotone lasso and general loss function, we have no exact algorithms for tracking the path, and hence for finding solutions at any point on the path. Friedman & Popescu (2004), however, have developed efficient $\epsilon$-stepping algorithms for finding forward-stagewise solutions for a variety of losses.

## 5. Discussion of Criteria

In conducting this research, we had several interesting false starts in terms of finding a criterion for forward stagewise. We briefly discuss some of these here.

We saw in the top right panel of Figure 3 on page 12 the residual sum of squares

$$\text{RSS}(\ell) = \sum_{i=1}^{N} \Big(y_i - \sum_{j=1}^{p} x_{ij}\beta_j(\ell)\Big)^2$$

for each of the three methods, as a function of their $L_1$ arc-length $\ell$. The curve for the forward stagewise sequence always lies below the curves for the other two methods. We were able to show a local optimality for forward stagewise in Theorem 2.

We initially had thought that forward stagewise might enjoy a global optimality criterion like the lasso.

**Candidate criterion 1:** *For each $L_1$ arc-length $\ell$, the forward stagewise coefficient $\beta(\ell)$ minimizes $\text{RSS}(\ell)$.*

This is true if the lasso paths are monotone, because then the two procedures coincide, as do $L_1$ arc-length and $L_1$ norm. But in general it is not the case.

**Lemma 2.** *In general there does not exist a coefficient profile that for all $\ell$ minimizes $\text{RSS}(\ell)$ over the set of curves having $L_1$ arc-length at most $\ell$.*

*Proof.* For any $\ell$ construct the "unit speed" coefficient path from the origin to the lasso solution for that $\ell$. This has $L_1$ arc-length and $L_1$ norm equal to $\ell$, and hence has the minimum value of $\mathrm{RSS}(\ell)$ over all curves having $L_1$ arc-length $\ell$. Thus any solution to our problem must agree with the lasso solution for all $\ell$. From the right-hand panel of Figure 3, this is not the case when the lasso and forward stagewise profiles are different. □

Another attempt at a global formulation of the problem involve used the integrated loss. Define the set of monotone-increasing functions

$$D_L^M \triangleq \left\{ \beta : [0, L] \to \mathbb{R}^{2p} \mid \beta(0) = 0, \mathrm{TV}(\beta, \ell) \leq \ell, \ \forall \ell \leq L, \ \beta_j(\ell) \textit{non-decreas.} \right\},$$

having arc-length at most $\ell$ up to the point $\ell$, for all $\ell \leq L$. Since the class is monotone, the arc-length is the $L_1$ norm, so we are asking for a monotone path that simultaneously solves a sequence of lasso problems, subject to that constraint.

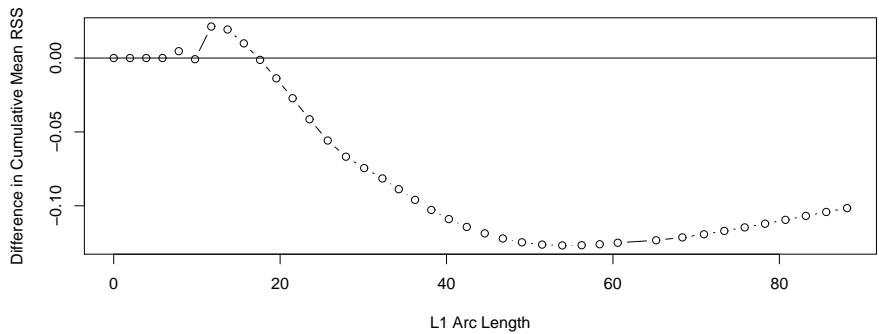

FIG 4. *A simulation provides a counter example to candidate 2. Shown is the difference between mean cumulative RSS for the exact solution to criterion 2 and forward stagewise (the former computed on a discretized set of 40 values for arc-length). Initially forward stagewise wins, only to be overtaken by the exact solution.*

**Candidate criterion 2:** *The forward stagewise algorithm minimizes the integrated residual sum of squares over the monotone class $D_L^M$:*

$$\tilde{\beta}(L) = \mathrm{argmin}_{\beta \in D_L^M} \int_0^L \sum_{i=1}^N \Big( y_i - \sum_{j=1}^p \tilde{x}_{ij} \beta_j(\ell) \Big)^2 d\ell \qquad (34)$$

As the integrated loss is a continuous, strictly convex functional on $D_L^M$, there exists a unique optimal path that solves (34). However it turns out that the forward stagewise solution is not always the optimal path.

We computed the exact solution to (34) for our simulation example, at a discretized sequence of 40 values for arc-length. Figure 4 compares the results with the forward stagewise solution at these same points. We compute for each the cumulative mean RSS, and plot their difference (exact−FS$_0$). In keeping with its greedy nature, FS$_0$ initially wins, only to be overtaken by the exact procedure. Hence FS$_0$ does not in general optimize criterion 2.

## 6. Monotonicity of Profiles

We have yet to say how the example of Figure 1 was generated. The data were generated from the model

$$Y = \sin(6X)/(1+X) + Z/4, \tag{35}$$

with $X$ taking 300 equally spaced values in $[0,1]$ and $Z \sim N(0,1)$. The 10 predictors were piecewise linear basis functions $(x-t_k)\cdot I(x > t_k)$, for each of the knots $\{t_k\}_1^{10} = \{0.0, 0.1, 0.2, \ldots 0.9\}$. Figure 5 shows the successive approximations to $\sin(6x)/(1+x)$ by the different methods, for five equally spaced solutions along their paths. Despite the differences in their coefficient profiles, the fits appear to be quite similar. The last column of Figure 5 uses piecewise constant basis functions $I(x > t_k)$ in place of the piecewise linear ones $(x - t_k) \cdot I(x > t_k)$. Figure 6 shows their coefficient profiles. Notice that all profiles are monotone, and hence the profiles for all three algorithms coincide.

The fact that they are the same under monotonicity is not a coincidence, and follows from their definitions. In that case there are no zero-crossing events, and then LAR and lasso coincide. In addition, monotonicity means that positive coefficients are never decreased and vice-versa, hence the non-negative least squares move in the forward stagewise procedure is the same as the least squares move in LAR.

Hence it is useful to characterize situations in which the coefficients profiles are monotone. Let $X$ denote the $N \times p$ matrix of standardized predictors, and let $X_A$ denote a subset of the columns of $X$, each multiplied by a set of arbitrary signs $s_1, s_2, \ldots s_{|A|}$. Finally, let $S_A$ be a diagonal matrix of the $s_j$ values. The results of Efron et al. (2004) show that a necessary and sufficient condition for every path to be monotone is

$$S_A(X_A^T X_A)^{-1} S_A \mathbf{1} \geq 0 \ \forall A \subseteq \{1, \ldots, p\}, \ S_A \tag{36}$$

In other words, for all subsets of predictors and sign changes of those predictors, the inverse covariance matrix must be diagonally dominant (this means that each diagonal element is at least as big as the sum of the of the other elements in its row).

For the piecewise-linear basis functions, it is easy to find a violation of (36); $A = \{4, 10, 9\}$ (the first three variables entered in Figure 1), with $S_A = \text{diag}(-1, 1, 1)$, gives some negative entries. Condition (36) clearly holds for any orthogonal basis, such as the Haar basis for piecewise constant fits. However,

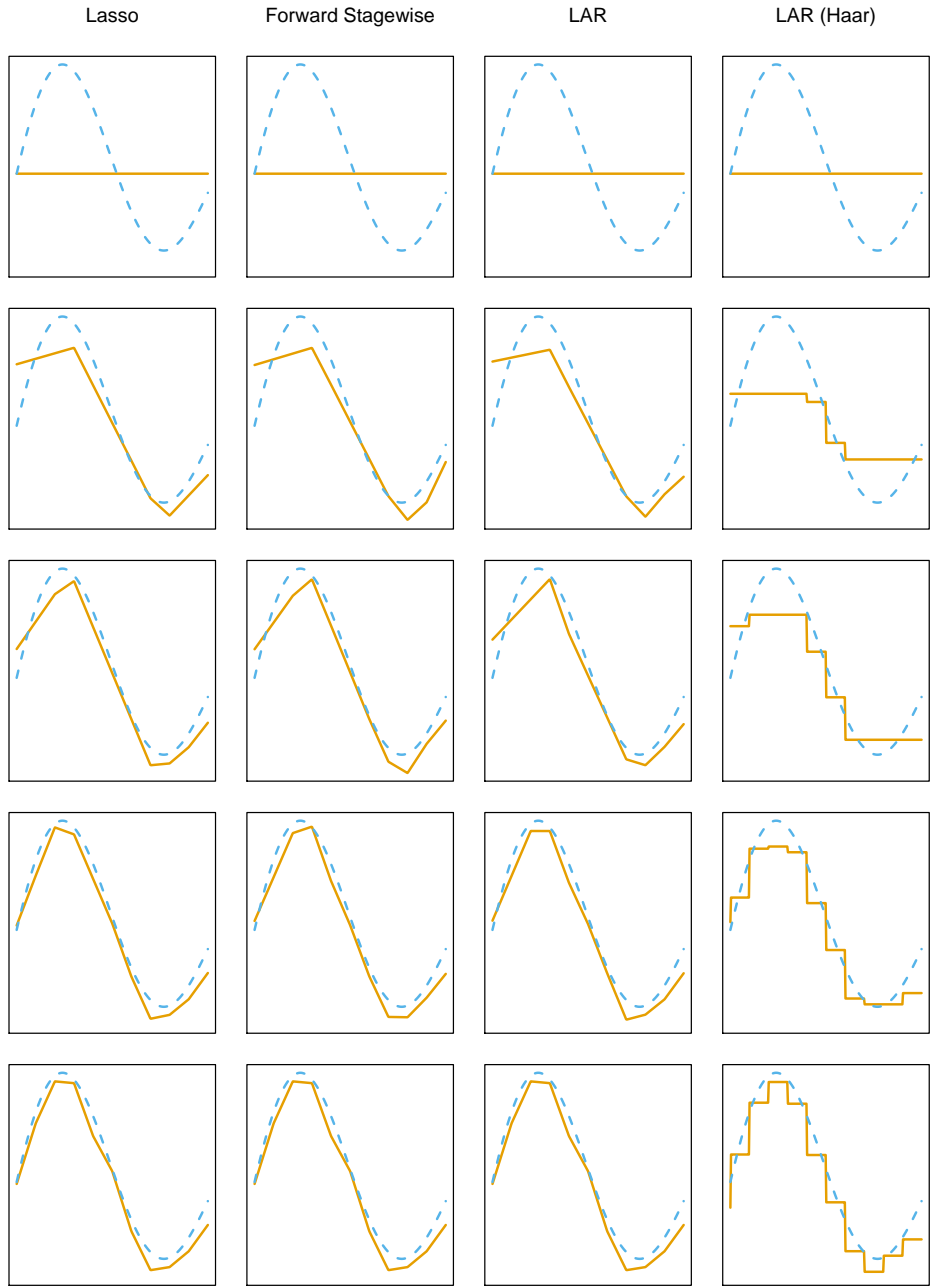

FIG 5. *Successive approximations to* $\sin(6x)/(1+x)$ *for five equally spaced solutions along their paths, for the example of Figure 1. The first three columns use piecewise linear bases; the last column uses piecewise constant bases, and the three methods coincide.*

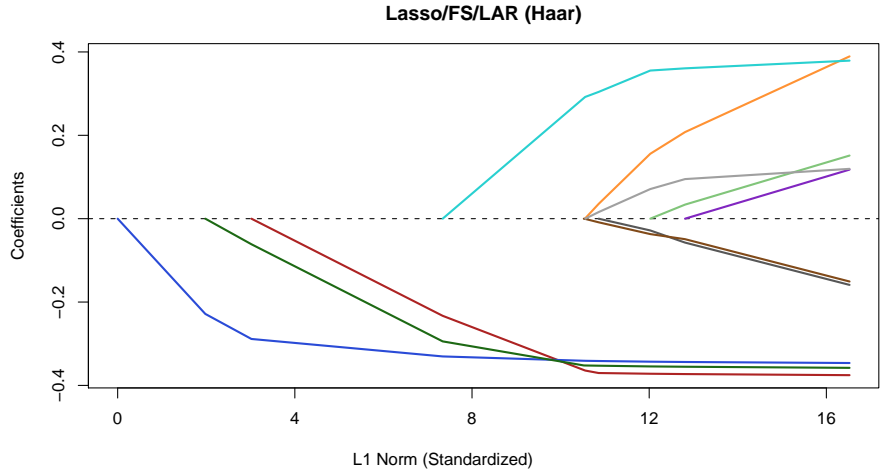

FIG 6. *Coefficient profiles for the same data as Figure 1, except that we have used piecewise constant basis functions. The coefficients are monotone, and lasso, $FS_0$ and LAR coincide.*

our piecewise constant basis is not orthogonal. We prove the following theorem in the appendix.

**Theorem 3.** *Condition (36) holds for piecewise constant bases, and hence the lasso, $FS_0$ and LAR solutions coincide.*

## 7. Lasso versus Forward Stagewise: Which is Better?

As discussed in Section 2, the current interest in $FS_\epsilon$ is because of its connection to least squares boosting. By understanding its properties in this simplified setting, we hope to learn more about the regularization path of boosting.

The results of this paper show that forward stagewise behaves like a monotone version of the lasso, and is locally optimal with regard to $L_1$ arc-length. This is in contrast to the lasso, which is less constrained.

This begs the question: with a large number of predictors, which algorithm is better? The monotone lasso will tend to slow down the search, not allowing the sudden changes of direction that can occur with the lasso. Is this a good thing?

To investigate this, we carried out a simulation study. The data consists of $N = 60$ observations on each of $p = 1000$ (Gaussian) variables, strongly correlated ($\rho = 0.95$) in groups of 20. The true model has nonzero coefficients for 50 variables, one drawn from each group, and the coefficient values themselves are drawn from a standard Gaussian. Finally Gaussian noise is added with variance $\sigma^2 = 36$, resulting in a noise-to-signal ratio of about 0.72. See Appendix A.3 on page 27 for more details.

The grouping of the variables is intended to mimic the correlations of nearby trees in boosting, and with the forward stagewise algorithm, this setup is intended as an idealized version of gradient boosting with shrinkage.

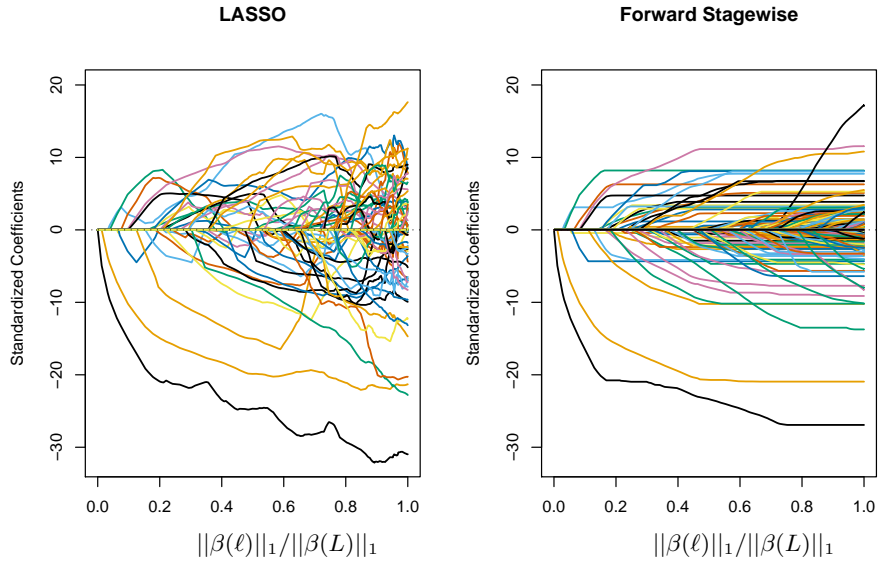

FIG 7. *Comparison of lasso and forward stagewise paths on simulated regression data. The number of samples is 60 and the number of variables 1000. The forward-stagewise paths fluctuate less than those of lasso in the final stages of the algorithms. Both paths are indexed by $L_1$-norm $||\beta(\ell)||_1$, scaled as a fraction of the $L_1$-norm at the end of the path $||\beta(L)||_1$.*

Figure 7 shows the coefficient paths for lasso and forward stagewise for a single realization from this model.

Here the coefficient profiles are similar only in the early stages of the paths. For the later stages, the forward stagewise paths are much smoother — in fact exactly monotone here — while those for the lasso fluctuate widely. This is due to the strong correlations among subsets of the variables.

The test-error performance of the two models is rather similar (figure 8), and they achieve about the same minimum. In the later stages forward stagewise takes longer to overfit, a likely consequence of the smoother paths. We are using $||\beta(\ell)||_1$ to index both curves for this plot, which would look quite different if instead we used arc-length. Since the forward stagewise path is monotone here, the $L_1$ norm is arc-length, so in a sense both MSE profiles are measured using their appropriate index.

On a more theoretical note, Buhlmann (2006) proves consistency of forward stagewise for high-dimensional linear models. See also Tropp (2004) and Tropp (2006) for comparisons of lasso and forward stagewise regression.

We conclude that for problems with large numbers of correlated predictors, the forward stagewise procedure and its associated $L_1$ arc-length criterion might be preferable to the lasso and $L_1$ norm criterion. This suggests that for general boosting-type applications, the incremental forward stagewise algorithms which are currently used, might be preferable to algorithms that try to solve the equivalent lasso problem.

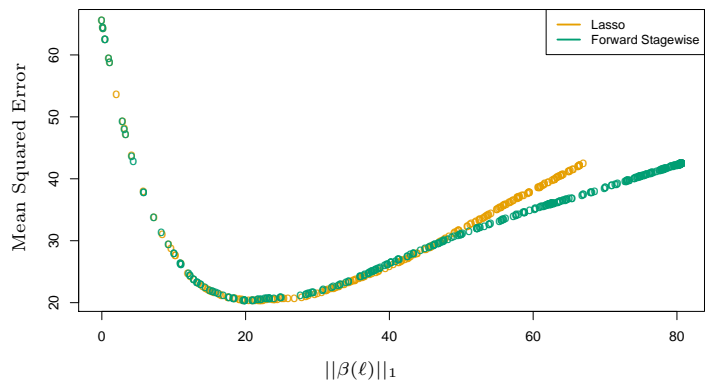

FIG 8. *Mean squared error for lasso and forward stagewise on the simulated data. Despite the difference in the coefficient paths, the two models perform similarly over the critical part of the regularization path. In the right tail, lasso appears to overfit more rapidly.*

## Appendix A: Appendix

### A.1. Proof of Theorem 1

*Part 1..* The Lagrangian corresponding to (3) is

$$\sum_i \left( y_i - \beta_0 - \left[ \sum_{j=1}^p x_{ij}\beta_j^+ - \sum_{j=1}^p x_{ij}\beta_j^- \right] \right)^2 + \lambda \sum_{j=1}^p (\beta_j^+ + \beta_j^-)$$
$$- \sum_{j=1}^p \lambda_j^+ \beta_j^+ - \sum_{j=1}^p \lambda_j^- \beta_j^-, \qquad (37)$$

with KKT conditions (for each $j$):

$$-\mathbf{x}_j^T \mathbf{r} + \lambda - \lambda_j^+ = 0, \qquad (38)$$
$$\mathbf{x}_j^T \mathbf{r} + \lambda - \lambda_j^- = 0, \qquad (39)$$
$$\lambda_j^+ \beta_j^+ = 0, \qquad (40)$$
$$\lambda_j^- \beta_j^- = 0. \qquad (41)$$

Here $\mathbf{r} = \mathbf{y} - \left[ \sum_{j=1}^p \mathbf{x}_j\beta_j^+ - \sum_{j=1}^p \mathbf{x}_j\beta_j^- \right]$ is the residual vector. From these we can deduce the following aspects of the solution:

1. If $\lambda = 0$, then $\mathbf{x}_j^T \mathbf{r} = 0 \forall j$, and the solution corresponds to the unrestricted least-squares fit.

2. 
$$\begin{aligned} \beta_j^+ > 0,\ \lambda > 0 &\implies \lambda_j^+ = 0 \\ &\implies \mathbf{x}_j^T \mathbf{r} = \lambda > 0 \\ &\implies \lambda_j^- > 0 \\ &\implies \beta_j^- = 0. \end{aligned}$$

3. Likewise $\beta_j^- > 0,\ \lambda > 0 \implies \beta_j^+ = 0$. Hence 2 and 3 give the intuitive result that only one of the pair $(\beta_j^+, \beta_j^-)$ can be positive at any time.

4. $|\mathbf{x}_j^T \mathbf{r}| \leq \lambda$.

5. If $\beta_j^+ > 0$, then $\mathbf{x}_j^T \mathbf{r} = \lambda$, or if $\beta_j^- > 0$, then $-\mathbf{x}_j^T \mathbf{r} = \lambda$

Since $\beta^0$ is on the lasso path, there exists a $\lambda = \lambda_0$ such that $\beta^0 = \beta(\lambda_0)$. Define the *active set* $\mathcal{A}$ to be the set of indices of variables in $\tilde{\mathbf{X}}$ with positive coefficients at $\lambda = \lambda_0$, and assume at $\lambda_1 = \lambda_0 - \Delta$ for suitably small $\Delta > 0$ this set has not changed. Define $\beta_\mathcal{A}(\lambda)$ to be the corresponding coefficients at $\lambda$. Then from deduction 5

$$\tilde{\mathbf{X}}_\mathcal{A}^T \left(\mathbf{y} - \tilde{\mathbf{X}}_\mathcal{A} \beta_\mathcal{A}(\lambda)\right) = \lambda \mathbf{1} \quad \text{for } \lambda \in [\lambda_1, \lambda_0]. \tag{42}$$

Hence

$$\tilde{\mathbf{X}}_\mathcal{A}^T \tilde{\mathbf{X}}_\mathcal{A} (\beta_\mathcal{A}(\lambda_1) - \beta_\mathcal{A}(\lambda_0)) = \Delta \mathbf{1}, \tag{43}$$

or

$$\beta_\mathcal{A}(\lambda_1) - \beta_\mathcal{A}(\lambda_0) = \Delta \cdot (\tilde{\mathbf{X}}_\mathcal{A}^T \tilde{\mathbf{X}}_\mathcal{A})^{-1} \mathbf{1}. \tag{44}$$

So while $\mathcal{A}$ remains constant, the coefficients $\beta_\mathcal{A}(\lambda)$ change linearly, according to (44). Since $\mathbf{r} = \mathbf{y} - \tilde{\mathbf{X}}_\mathcal{A} \beta_\mathcal{A}(\lambda_0)$ and $\tilde{\mathbf{X}}_\mathcal{A}^T \mathbf{r} = \lambda_0 \mathbf{1}$ from (42),

$$\beta_\mathcal{A}(\lambda_1) - \beta_\mathcal{A}(\lambda_0) = \frac{\Delta}{\lambda_0} \cdot (\tilde{\mathbf{X}}_\mathcal{A}^T \tilde{\mathbf{X}}_\mathcal{A})^{-1} \tilde{\mathbf{X}}_\mathcal{A}^T \mathbf{r} \tag{45}$$

as claimed.

The active set $\mathcal{A}$ will change if a variable "catches up" (in terms of 5), in which case it is augmented and the direction (44) is recalculated. It will also change if a coefficient attempts to pass through zero, in which case it is removed from $\mathcal{A}$. □

*Part 2..* At each step, the monotone incremental forward stagewise algorithm (Algorithm 3) selects the variable having largest correlation with the residuals, and moves its coefficient up by $\epsilon$. There may be a set $\mathcal{A}$ of variables competing for this maximal correlation, and a succession of $N$ such moves can be divided up according to the $N_j = \rho_j N$ that augmented variable $j$s coefficient. In the limit as $\epsilon$ decreases and $N$ increases such that $N\epsilon = \varepsilon$, we can expect an active set $\mathcal{A}$ of variables *tied* in terms of the largest correlation, and a sequence of moves of total $L_1$ arc-length $\varepsilon$ distributed among this set with proportions $\rho_\mathcal{A}$. Efron et al. (2004) showed in their Lemma 11 that for sufficiently small $\varepsilon$, this set would not change. Suppose $\tilde{\mathbf{X}}_\mathcal{A}^T \mathbf{r} = c\mathbf{1}$ for some $c$, reflecting the equal correlations.

The limiting sequence of moves $\varepsilon \rho_\mathcal{A}$ must have positive components, must maintain the equal correlation with the residuals, and subject to these constraints should decrease the residual sum-of-squares as fast as possible.

T. Hastie et al./The monotone lasso    25Consider the optimization problem

$$\min_\rho \frac{1}{2}||\mathbf{r} - \varepsilon\tilde{\mathbf{X}}_\mathcal{A}\rho||_2^2 \text{ subject to } \rho_j \geq 0, \sum_{j \in \mathcal{A}} \rho_j = 1, \quad (46)$$

which captures the first and third of these requirements. The Lagrangian is

$$L(\rho, \gamma, \lambda) = \frac{1}{2}||\mathbf{r} - \varepsilon\tilde{\mathbf{X}}\rho||_2^2 - \sum_{j=1}^{p} \gamma_j \rho_j + \lambda(\sum_j \rho_j - 1), \quad (47)$$

with KKT conditions

$$-\varepsilon\tilde{\mathbf{X}}_\mathcal{A}^T(\mathbf{r} - \varepsilon\tilde{\mathbf{X}}_\mathcal{A}\rho) - \gamma + \lambda\mathbf{1} = 0 \quad (48)$$

$$\begin{aligned} \gamma_j &\geq 0 \\ \rho_j &\geq 0 \\ \gamma_j \rho_j &= 0 \\ \sum_j \rho_j &= 1. \end{aligned} \quad (49)$$

Here $\gamma$ is a vector with components $\gamma_j$, $j \in \mathcal{A}$. Note that for $\rho_j > 0$, $\gamma_j = 0$, and hence (48) shows that the correlations with the residual remain equal, which was the second requirement above.

Consider a second optimization problem (the one in the statement of the the theorem):

$$\min_\theta \frac{1}{2}||\mathbf{r} - \tilde{\mathbf{X}}_\mathcal{A}\theta||_2^2 \text{ subject to } \theta_j \geq 0. \quad (50)$$

The corresponding KKT conditions are

$$-\tilde{\mathbf{X}}_\mathcal{A}^T(\mathbf{r} - \tilde{\mathbf{X}}_\mathcal{A}\theta) - \nu = 0 \quad (51)$$

$$\begin{aligned} \nu_j &\geq 0 \\ \theta_j &\geq 0 \\ \nu_j \theta_j &= 0 \end{aligned} \quad (52)$$

We now show that the $\hat{\rho} = \hat{\theta}/||\hat{\theta}||_1$ solves (46) for all $\varepsilon > 0$, where $\hat{\theta}$ is the solution to (50). From (48) we get

$$\varepsilon^2 \tilde{\mathbf{X}}_\mathcal{A}^T \tilde{\mathbf{X}}_\mathcal{A} \rho = (\lambda + \varepsilon c)\mathbf{1} + \gamma, \quad (53)$$

and from (51)

$$\tilde{\mathbf{X}}_\mathcal{A}^T \tilde{\mathbf{X}}_\mathcal{A} \hat{\theta} = c\mathbf{1} + \nu. \quad (54)$$

With $s = ||\hat{\theta}||_1 = \sum_{j \in \mathcal{A}} \hat{\theta}_j$, we multiply (54) by $\varepsilon^2/s$ to get

$$\varepsilon^2 \tilde{\mathbf{X}}_\mathcal{A}^T \tilde{\mathbf{X}}_\mathcal{A} \hat{\rho} = \frac{\varepsilon^2}{s}(c\mathbf{1} + \nu) \quad (55)$$
$$= (\lambda^* + \varepsilon c)\mathbf{1} + \gamma^*, \quad (56)$$

where

$$\gamma_j^* = \nu_j \frac{\varepsilon^2}{s}$$
$$\lambda^* = \frac{\varepsilon^2 c}{s} - \varepsilon c.$$

It is easy to check, using (52), that $(\hat{\rho}, \gamma^*, \lambda^*)$ satisfy (48)-(49).

Variables with $\hat{\theta}_j = 0$ may drop out of the active set, since from (48) and (49), if $\gamma_j > 0$, for $\varepsilon > 0$ their correlation will decrease faster than those with positive coefficients.

This directions is pursued until a variable not in $\mathcal{A}$ "catches up" in terms of correlation, at which point the procedure stops, $\mathcal{A}$ is updated, and the direction is recomputed. □

### A.2. Proof of Theorem 3

We need to verify that when using piecewise constant basis functions, (36) holds for every $A$ and every sign matrix $S_A$.

Suppose we use $k$ piecewise constant basis functions with knots $t_1 \leq \cdots \leq t_k$. Let

$$n_j = \sum_{i=1}^{n} I(x_i > t_j)$$

be the number of observed $x$'s to the right of the $j$-th knot. Without loss of generality, we assume that each $n_j > 0$, otherwise that predictor contributes nothing to the model as all observed $x$'s are either to the right or the left of that knot.

A simple calculation shows that, for $i \leq j$, after normalizing the columns of $X$,

$$(X^t X)_{ij} = \sqrt{\frac{(n - n_i)}{n_i} \cdot \frac{n_j}{n_j - n}},$$

which is the covariance function of a Brownian bridge $(B_s)_{0 \leq s \leq 1}$, normalized to have unit variance, at the time points $0 < s_1 \leq \cdots \leq s_k < 1$

$$s_j = \frac{n_{k-j+1}}{n}.$$

If we can prove that $(X^t X)^{-1}$ is diagonally dominant for every $k$ and every choice of knots, then, as every principal minor is of exactly the same form (with a smaller $k$ and fewer knots), we will also have proved that $(X_A^t X_A)^{-1}$ is diagonally dominant, hence that the lasso paths are monotone.

We prove that $(X^t X)^{-1}$ is diagonally dominant by computing $(X^t X)^{-1}$. One way to compute the $(X^t X)^{-1}$ is to compute the density of

$$\left( \frac{B_{s_1}}{\sqrt{s_1(1-s_1)}}, \ldots, \frac{B_{s_k}}{\sqrt{s_k(1-s_k)}} \right)$$

and read off the inverse from the exponent of the density. Before we turn to computing the density, we note that

$$\left( \frac{B_{s_1}}{\sqrt{s_1(1-s_1)}}, \ldots, \frac{B_{s_k}}{\sqrt{s_k(1-s_k)}} \right) \stackrel{D}{=} \left( \frac{W_{v_1}}{\sqrt{v_1}}, \ldots, \frac{W_{s_k}}{\sqrt{v_k}} \right)$$

where $W$ is a standard Brownian motion and

$$v_j = \frac{s_j}{1 - s_j}.$$

It is now simple to show that, up to a constant multiple, the exponent of the density of
$$\left(\frac{W_{v_1}}{\sqrt{v_1}}, \ldots, \frac{W_{v_k}}{\sqrt{v_k}}\right)$$
evaluated at $(w_1, \ldots, w_k)$ is
$$w_{s_1}^2 + \sum_{j=2}^{k} \frac{(v_j^{1/2} w_{v_j} - v_{j-1}^{1/2} w_{v_{j-1}})^2}{v_j - v_{j-1}}.$$

Therefore $(X^t X)^{-1}$ has elements
$$(X^t X)^{-1}_{j,j} = v_j \left(\frac{1}{v_j - v_{j-1}} + \frac{1}{v_{j+1} - v_j}\right)$$
$$= \frac{(v_{j+1} - v_{j-1}) v_j}{(v_j - v_{j-1})(v_{j+1} - v_j)}$$
$$(X^t X)^{-1}_{j,j-1} = -\frac{\sqrt{v_j v_{j-1}}}{v_j - v_{j-1}}.$$

Because the off-diagonal entries of $(X^t X)^{-1}$ are non-positive, we only have to show that
$$(X^t X)^{-1} \mathbf{1} \geq 0,$$
as multiplying on the left and right by $S_A$ will only increase the entries of the above vector.

We must therefore prove that for all $j$,
$$\frac{(v_{j+1} - v_{j-1}) v_j}{(v_j - v_{j-1})(v_{j+1} - v_j)} - \frac{\sqrt{v_j v_{j-1}}}{v_j - v_{j-1}} - \frac{\sqrt{v_{j+1} v_j}}{v_{j+1} - v_j} \geq 0.$$

By scaling and combining fractions, the above is implied by the following: for every $a < 1 < b$
$$1 - \sqrt{a} \cdot \frac{b-1}{b-a} - \sqrt{b} \cdot \frac{1-a}{b-a} \geq 0.$$

It remains therefore to prove that this inequality holds. However, this is just Jensen's inequality: define a two-point distribution placing mass $(b-1)/(b-a)$ on $\sqrt{a}$ and mass $(1-a)/(b-a)$ on $\sqrt{b}$. Then, if $Z$ is distributed according to this law:
$$\mathrm{E}(X) = \sqrt{a} \cdot \frac{b-1}{b-a} + \sqrt{b} \cdot \frac{1-a}{b-a} \leq (\mathrm{E}(X^2))^{1/2} = 1.$$

We note that general conditions for monotonicity can be derived. However it is not clear how these might be verified in practice.

### A.3. Simulation Details

Here we give more details of the simulation in Section 7. The regression model has the form
$$Y = X\beta + \epsilon, \tag{57}$$

where $X \sim N(0, \Sigma)$. $X$ has 1000 components, correlated in blocks of size 20. Hence $\Sigma$ is a block-diagonal covariance matrix, with 50 blocks $\Sigma_m$ each of size 20, and each block has the (identical) form

$$\Sigma_m = (1-\rho)\mathbf{I}_{20} + \rho\mathbf{1}\mathbf{1}^T. \tag{58}$$

We used $\rho = 0.95$ in each block. The 1000-vector $\beta$ was chosen to be sparse, with only 50 non-zero entries — one per block. Without loss of generality, we picked the first variable in each block to have a non-zero coefficient, drawn at random from a standard Gaussian distribution. The noise term $\epsilon$ was also Gaussian, with variance $\sigma^2 = 36$. $N = 60$ realizations were drawn from this model. For the noise-to signal ratio we compute $\text{var}(\epsilon)/\text{var}(X\beta)$. Since $X$ and $\beta$ have mean zero, the denominator is $E(\beta^T \Sigma \beta) = \text{tr}[\Sigma E(\beta\beta^T)] = 50$, hence the ratio is 0.72.